\documentclass[12pt,twoside]{article}
\usepackage{amsmath}

\date{}

\begin{document}

\emph{Journal Name, Vol. X, No. X, DATE YEAR, pp. xx-xx}

%\emph{ISSN 2219-7184; Copyright \copyright ICSRS Publication, 2010}

%\emph{www.i-csrs.org}

%\emph{Available free online at http://www.geman.in}

\centerline{}

\centerline{}

\centerline {\Large{\bf A Note on Generating Almost Pythagorean Triples}}

\centerline{}

\centerline{}

%% My definition
\newcommand{\mvec}[1]{\mbox{\bfseries\itshape #1}}

\centerline{\bf {John Rafael M. Antalan$^*$ and Mark D. Tomenes}}

\centerline{}

\centerline{Department of Mathematics and Physics} 
\centerline{College of Arts and Sciences}
\centerline{ Central Luzon State University}
\centerline{Science City of Munoz, Nueva Ecija, Philippines }

\centerline{E-mail $^*$: jrantalan@clsu.edu.ph}

\centerline{}

\newtheorem{Theorem}{\quad Theorem}[section]

\newtheorem{Definition}[Theorem]{\quad Definition}

\newtheorem{Corollary}[Theorem]{\quad Corollary}

\newtheorem{Lemma}[Theorem]{\quad Lemma}

\newtheorem{Example}[Theorem]{\quad Example}

\centerline{\bf Abstract}
{\emph{In 1987, Orrin Frink introduced the concept of almost Pythagorean triples. He defined them as an ordered triple $(x,y,z)$ that satisfies the Diophantine equation $x^2+y^2=z^2+1$ where $x,y$ and $z$ are positive integers. In his paper, he determined all almost Pythagorean triples by solving the given Diophantine equation. However, the result does not explicitly and readily give a particular almost Pythagorean triple. In this note,  using basic algebraic operations and Frink's result we give an explicit formula that readily gives a particular almost Pythagorean triple. We also give a certain integer sequence as a result of the generated formula.}}
\\


{\bf 2000 MSC No:} 11A99, 11D45.

%=============================
\section{Introduction}
%=============================

A Pythagorean triple (PT) is an ordered triple $(a,b,c)$ that satisfies the Pythagorean equation:
\begin{equation}
a^2+b^2=c^2.	
\end{equation}

If $a$,$b$ and $c$ are pairwise relatively prime then the triple is a primitive Pythagorean triple (PPT).
Modifying (1) slightly we have:
\begin{equation}
x^2+y^2=z^2+1.
\end{equation}

Let us agree to name (2) as ``almost Pythagorean equation". An ordered triple that satisfies (2) was named almost Pythagorean triple (APT) by O. Frink in [1]. He also proved that:

\begin{Theorem}
If $(a,b,c)$ is a primitive Pythagorean triple and if $(p,q,r)$ is an almost Pythagorean triple, then the triples:

\begin{equation}
(x,y,z)=(at+p,bt+q,ct+r)
\end{equation}

and

\begin{equation}
(x',y',z')=(at+p',bt+q',ct+r')
\end{equation}
are almost Pythagorean triples for all positive integer $t$ and for unique integers $p,p',q,q',r$ and $r'$ that depends on the primitive Pythagorean triple $(a,b,c)$ satisfying $p+p'=a$, $q+q'=b$ and $r+r'=c$. Likewise all almost Pythagorean triples takes the form (3) and (4).
\end{Theorem}
This result suggests that in order to generate APTs: 
\begin{enumerate}
\item First it is necessary to have a PPT $(a,b,c)$. 
\item Use the PPT, to generate the simultaneous system of 6 linear equations (with $t=1$):
\begin{center}
\[ \begin{cases}
p^2+q^2=r^2+1\\
(a+p)^2+(b+q)^2=(c+r)^2+1\\
(a+p')^2+(b+q')^2=(c+r')^2+1\\
p+p'=a\\
q+q'=b\\
 r+r'=c
\end{cases}\]
\end{center}
\item Lastly, solve the system in step 2 to find the integers $p,p',q,q',r$ and $r'$. 
\item Since the variables $a,b,c,p,p',q,q',r$ and $r'$ are now known constants, the  triples of the form (3) and (4) are now APT's. 
\end{enumerate}
For example, in order to generate the almost Pythagorean triples:

\begin{equation}
(x,y,z)=(3t+1,4t+3,5t+3)
\end{equation}
and 

\begin{equation}
(x',y',z')=(3t+2,4t+1,5t+2)
\end{equation}
\begin{enumerate}
\item The PPT $(a,b,c)=(3,4,5)$ was being considered. 
\item This PPT yields the simultaneous system of equations:
\begin{center}
\[ \begin{cases}
p^2+q^2=r^2+1\\
(3+p)^2+(4+q)^2=(5+r)^2+1\\
(3+p')^2+(4+q')^2=(5+r')^2+1\\
p+p'=3\\
q+q'=4\\
 r+r'=5
\end{cases}\]
\end{center}
\item The solutions of the last system of linear equations, can easily be solved and is given by $p=1$, $q=3$, $r=3$, $p'=2$, $q'=1$ and $r'=2$. 
\item By using the solved constants in (3) and (4) we are done.
\end{enumerate}

Based from the enumerated and illustrated steps however, the general solution of (2) does not explicitly and readily give a particular almost Pythagorean triple. Also for PPT (a,b,c) with large components, finding integers $p,p',q,q',r$ and $r'$ seems to be a not an easy task.
\\

So, in this note, using a certain formula that generates ``all" PPT [2], suggested steps in generating APTs from Theorem 1.1 and simple algebraic operations, we address the limitations being stated. That is, we give an explicit formula that generates APTs and show that this formula is ideal in generating APTs even if the components of the PPTs in which they were derived are large. 
\\

%==============================
\section{Result}
%==============================

Before presenting our result, we first state a preliminary lemma from [3]. 

\begin{Lemma}
All of the solutions of the Pythagorean equation $a^2+b^2=c^2$ satisfying the conditions:
\begin{enumerate}
\item gcd$(a,b,c)=1$
\item $2|b$
\item $a,b,c>0$
\end{enumerate}
is given by $(a,b,c)=(s^2-k^2, 2sk, s^2+k^2)$. For relatively prime integers $s>k>0$ and $s\neq k (mod\ 2)$.
\end{Lemma}

Now, let $i\geq 2$ be an integer. Note that $s_i=i$ and $t_i=i-1$  are relatively prime for all $i$ and it satisfies the incongruence relation $s_i\neq t_i(mod\ 2)$. Using lemma 2.1, we generate the PPT of the form $(2i-1,2i^2-2i,2i^2-2i+1)$. 
\\

Using the suggested method of Theorem 1.1 we have:  
\begin{enumerate}
\item Use the  PPT of the form $(2i-1,2i^2-2i,2i^2-2i+1)$ with $i\geq 2$.
\item We then generate the system of linear equations
\begin{center}
 \[\begin{cases}
p^2+q^2=r^2+1\\
[(2i-1)+p]^2+[(2i^2-2i)+q]^2=[(2i^2-2i+1)+r]^2+1\\
[(2i-1)+p']^2+[(2i^2-2i)+q']^2=[(2i^2-2i+1)+r']^2+1\\
p+p'=2i-1\\
q+q'=2i^2-2i\\
r+r'=2i^2-2i+1
\end{cases}
\] 
\end{center} 
\item Solving for the integers $p,p',q,q',r$ and $r'$ and by letting $(p,q,r)$ to be APT we have the following computations:
\begin{center}
$a=2i-1$; $b=2i^2-2i$ and $c=2i^2-2i+1$
\end{center}
Using the second equation above we have:
\begin{center}
$(a+p)^2+(b+q)^2=(c+r)^2+1$
\end{center}
Expanding the last displayed equation will lead to:
\begin{center}
$a^2+2ap+p^2+b^2+2bq+q^2=c^2+2cr+r^2+1$
\end{center}
Using the first equation above and by noting that $(a,b,c)$ is a PPT, the last displayed equation becomes:
  \begin{center}
$2ap+2bq=2cr$
\end{center}
After some simplifications and expressing $a,b$ and $c$ in terms of $i$, we have:
\begin{center}
$(2i-1)p+(2i^2-2i)q=(2i^2-2i+1)r$
\end{center}
Rewriting then and simplyfying we have:
\begin{center}
$(2i^2-2i)q+(2i-1)p=(2i^2-2i)r+r$
\end{center}
Finally, we get the clear solutions $p=1$, $q=2i-1$ and $r=2i-1$. And using the fourth, $5^{th}$ and $6^{th}$ equations we conclude that $p'=2i-2$, $q'=2i^2-4i+1$ and $r'=2i^2-4i+2$.
\end{enumerate}

%---------------------------------------------------------------------------------------------------------------------------------------------------
Thus we proved the theorem:
\begin{Theorem}
The triples given by:
\begin{equation}
(x,y,z)=(at+1,bt+(2i-1),ct+(2i-1))
\end{equation}
and
\begin{equation}
(x',y',z')=(at+(2i-2),bt+(2i^2-4i+1),ct+(2i^2-4i+2))
\end{equation}\\
where\begin{center}$(a,b,c)=(2i-1,2i^2-2i,2i^2-2i+1)$ ,$i\geq 2$,\end{center} are almost Pythagorean triples for all $t\in Z^+$. 
\end{Theorem}

%=============================================
\section{Examples}
%=============================================

In this section, we illustrate how to use Theorem 2.2 in generating almost Pythagorean triples. Recall that the  formulas:
\begin{center}
$(x,y,z)=(at+1,bt+(2i-1),ct+(2i-1))$\\
$(x',y',z')=(at+(2i-2),bt+(2i^2-4i+1),ct+(2i^2-4i+2))$
\end{center}
generate APTs with $(a,b,c)=(2i-1,2i^2-2i,2i^2-2i+1)$ for $i\geq 2$ and $t$ a positive integer. Let us use this formula to generate some APTs.  
\\

{\bf Example 1:} Let $i=4$ and $t_1=5$ and $t_2=6$. 
\\

When $i=4$ and $t=5$, the formula above generates the almost Pythagorean triples: 
\begin{center}
$(36,127,132)$

$(41,137,143)$. 
\end{center}
When $t=6$ we generate the APT's 
\begin{center}
$(43,151,157)$\\
$(48,161,168)$
\end{center}

{\bf Example 2:} Let $i_1=10$, $i_2=11$ and $t=7$.
\\

When $i=10$ and $t=7$ the formulas above yields the APTs 
\begin{center}
$(134,1279,1286)$\\
$(151,1421,1429)$.
\end{center}

 When $i=11$ and $t=7$ we have the APTs 
\begin{center}
$(148,1561,1568)$\\
$(167,1739,1747)$.
\end{center}

{\bf Example 3:} Let $i=3120$ and $k=25$. 
\\

Letting $i=3120$ and $t=25$, the formulas above yields the triples:  \begin{center}
$(155\ 976, 486\ 570\ 239, 486\ 570\ 264)$
\end{center}
and
\begin{center}
$(162\ 213, 506\ 020\ 321, 506\ 020\ 347)$.
\end{center}
and they are almost Pythagorean triples which can easily be verified using any computing software.

%=============================================
\section{A Certain Integer Sequence}
%=============================================
Comparing the almost Pythagorean triples generated by (7) and (8) we see that:
\begin{center}
$(x,y,z)<(x',y',z')$
\end{center}
setting $t=1$ and running $i$ over the set of integers greater than or equal to two  we have an infinite almost Pythagorean triples the first 10 of which are shown on the table below:
\begin{center}
\begin{tabular}{||c c c||}
\hline 
x & y & z
\\ [0.5ex]
\hline \hline
4 & 7 & 8\\
6 & 17 & 18\\
8 & 31 & 32\\
10 & 49 & 50\\
12 & 71 & 72\\
14 & 97 & 98\\
16 & 127 & 128\\
18 & 161 & 162\\
20& 199 & 200\\
22 & 241 & 242\\
[1ex]
 \hline
\end{tabular}
\end{center}
Rewriting the triples horizontally we have the sequence of integers:
\begin{center}
4, 7, 8, 6, 17, 18, 8, 31, 32, 10, 49, 50, 12, 71, 72, ...
\end{center}
prior to the creation of this version of this manuscript, the sequence generated above was applied as new integer sequence on The On-line Encyclopedia of Integer Sequences (OEIS) and was accepted and given the sequence number {\bf A261654} [4] and is called ``Lead almost-pythagorean triples generated by primitive pythagorean triples of the form $(2i-1,2i^2-2i,2i^2-2i+1), i\geq 2$".     
%=============================================
\section{Concluding Remark}
%=============================================
In this note, we successfuly gave an explicit formula in generating almost Pythagorean triples. These formulas were stated in Theorem 2.2 and in the beginning paragraph of Example section. The formulas suggest that in order to generate APTs just think of any positive integer $t$ and an integer $i\geq 2$, after using them in the formulas the output is the desired APT.  
\\

However, the result in Theorem 2.2 does not generate all almost Pythagorean triples explicitly since we restrict there our PPTs to be of the form $(a,b,c)=(2i-1,2i^2-2i,2i^2-2i+1)$.\\

As a possible extension to this note, we encourage the readers to answer the problem: ``Is there an explicit formula that generates all almost Pythagorean triples? If none, given two explicit solutions (comes in pair) $s_1$ and $s_2$ , can we possibly device a criterion that determines which one is ideal, perhaps in terms of the number of APT's being generated?". 
\bigskip

{\bf ACKNOWLEDGEMENTS}
\\

The authors highly acknowledge the help of their colleagues, friends and the Central Luzon State University in general for their encouragement and support. The first author is highly indebted to his love Josephine Joy Tolentino for reading the manuscript and suggesting some changes in the format of the paper and the words being used.                 

\end{document}